\input amstex 
\documentstyle{amsppt}
\input bull-ppt
\keyedby{bull308e/mhm}

\define\xt{\operatorname{xt}}
\define\diam{\operatorname{diam}}
\define\curv{\operatorname{curv}}
\define\pack{\operatorname{pack}}
\define\Sec{\operatorname{sec}}
\define\rad{\operatorname{rad}}
\define\exc{\operatorname{exc}}
\define\dist{\operatorname{dist}}

\topmatter
\cvol{27}
\cvolyear{1992}
\cmonth{October}
\cyear{1992}
\cvolno{2}
\cpgs{261-265}
\title Curvature, triameter, and beyond \endtitle
\author Karsten Grove and Steen Markvorsen \endauthor
\address Department of Mathematics,
University of Maryland, College Park, Maryland 
20742-0001\endaddress
\ml kng\@math.umd.edu \endml
\address Mathematics Institute, 
Technical University of Denmark, 2800 Lyngby, 
Denmark\endaddress
\ml steen\@mat.dth.dk \endml
\date November 5, 1991\enddate
\subjclassrev{Primary 53C20, 51K10, 53C23}
\thanks Both authors were supported in part 
by the Danish Research Council; 
the first author was also supported 
by a grant from the National Science Foundation\endthanks
\abstract In its most general form, the recognition 
problem in Riemannian
geometry asks for the identification of an unknown 
Riemannian manifold via
measurements of metric  invariants on the manifold.  We 
introduce a new
infinite sequence of invariants, the first term of which 
is the usual diameter,
and illustrate the role of these global shape invariants 
in a number of
recognition problems. \endabstract
\endtopmatter

\document
	
It is apparent that information about basic geometric 
invariants such as  {\it
curvature, diameter\/}, and {\it volume\/} alone does not 
suffice in order to
characterize Riemannian manifolds in general.  For this 
reason it is not only
natural but necessary to pursue and investigate other 
metric invariants. 
Recently such investigations have included Gromov's {\it 
filling radius\/}
(cf., e.g., \cite {G, K1, K2,  W}), {\it excess\/} 
invariants
(cf., e.g., \cite {GP2, O, PZ}), and Urysohn's  {\it
intermediate diameters\/} (cf., e.g., \cite {U, G, K3}).  

	The purpose here is to introduce new metric invariants 
and announce
related {\it recognition theorems}.  Our main concern is 
an infinite sequence: 
{\it diameter, triameter, quadrameter, 
quintameter\/},$\dotsc$, etc.,
whose $n$\<th term is
based on measurements on $(n+1)$-tuples of points.  
Precisely, if $(X,\dist)$ is
any compact metric space, the $q$-{\it extent\/}, 
$\xt_qX$, of $X$ is the maximal
average distance between $q$ points in $X$, i.e.,
$$ 
	\text{xt}_qX = \max_{(x_1,\dotsc,x_q)} \sum_{i<j} 
\text{dist}(x_i,x_j)/
\binom q2, \qquad \;  x_1,\dotsc,x_q \in X.
$$
With this definition, $\xt_2X = \text{diam} X, 
\text{xt}_3X = \text{triam}X$,
etc., and obviously
$$
	\text{xt}_2X \geq \text{xt}_3X \ge\cdots\geq \text{xt}_qX 
\geq
\text{xt}_{q+1}X \geq\cdots\geq \text{xt}X,
$$
where $\xt X = \lim_q \; \text{xt}_qX$ 
is called the {\it extent\/}  of $X$.  It is easy to see 
that
$$
\tsize\frac 12 \text{diam} X \le \text{xt}X < \text{diam}X
$$
for any compact $X$ and that these inequalities are 
optimal.  Somewhat 
surprisingly, however, it turns out 
that $\xt X$ is related to the excess, $\exc
X$,  as defined in \cite {GP2}.  Namely, if $X$ has almost 
minimal extent, its
excess is almost zero.  From this and several explicit 
computations (cf.\ \cite
{GM}), it appears that the above extent invariants
are sensitive to asymmetries of a space and, therefore, 
should be thought of
not as size invariants, but as {\it global shape 
invariants\/} in the same way
that curvatures are thought of a {\it local shape 
invariants\/}.

\heading Spaces with large extents \endheading
The first step in the recognition program for a given 
metric invariant 
is to investigate its range when restricted to various 
subclasses of metric
spaces  (cf.\ \cite {G2, GM}).  For the individual 
$q$-extents there are
trivial optimal estimates like the one given for $\xt X$ 
above, even when
restricted to Riemannian manifolds.   This becomes an 
attractive
problem, however, when a lower curvature bound is present, 
and thus local and global
shape invariants are balanced against one another.

Specifically, for each $k \in \Bbb R$ and $R > 0$, let 
$D_k^n(R)$ denote the
closed metric $R$-ball in the simply connected 
$n$-dimensional complete space
form $S_k^n$ of constant curvature $k$.  For any closed 
Riemannian
$n$-manifold $M$, whose {\it sectional curvature\/} and 
{\it radius\/} satisfy
$\Sec M \ge k$ and $\rad M \le R$, standard Toponogov 
distance comparison
yields 
$$
\text{xt}_q \; M \le \text{xt}_q \; D_k^n(R),
\tag {$*$}
$$
for any integer $q \ge 2$.  Recall that $\rad X \le R$ 
if and only if $X = D(x,R)$ for some $x \in X$ (cf.\ \cite 
{SY, GP4}).

When $k > 0$  and $R \ge \pi/2\sqrt k$, it turns out that 
$\xt_q
D_k^n(R) = \text{xt}_q S_k^n = \text{xt}_q[0,\pi/\sqrt k]$ 
for all $q$, and all
$q$-{\it extenders\/}, i.e., $q$-tuples of points 
realizing $\xt_q$, can be
explicitly described (cf.\ \cite {GM, N}). Moreover the 
inequalities $(*)$ 
are optimal in this case.  The following range/recognition 
theorem generalizes
Toponogov's maximal diameter theorem (cf.\ \cite {CE}):   

\thm{Theorem A}  Let $M$ be a closed Riemannian 
$n$-manifold, $n \ge 2$, with 
$\Sec M \ge 1$. 
	 Then 
$	\xt_q M \le \xt_q S_1^n$
	for every $q \ge 2$.

	If equality holds for some $q \ge 2$, then $M$ is 
isometric to $S_1^n$.

For any $\varepsilon > 0$ and any $\pi/2 \le R \le \pi$ 
there is a
Riemannian manifold $M \simeq S^n$ with $\Sec M \ge 1, 
\rad M \le R$, and
$\xt_q M \ge \xt_q S_1^n - \varepsilon$. \ethm

Note that since $\diam M \ge \text{xt}_q M$ for any $q \ge 
2$, the diameter
sphere theorem \cite {GS} implies that $M$ is homeomorphic 
to $S^n$ if  $\xt_q
M > \pi/2$ in the above theorem.

In the remaining cases we only know $\xt _q D_k^n(R)$ when 
$q \le n+1$ (cf.\
\cite {GM, T, H}).  Here $\xt _{n+1} D_k^n(R)$ is of
particular interest because for  $R < \pi/4 \sqrt k$ (if 
$k > 0)$ there is only
one $(n+1)$-extender, namely, the vertices of the unique 
maximal, regularly
inscribed $n$-simplex, $\Delta_k^n(R)$ in $D_k^n(R)$.  
Here without the loss of
generality we  may assume $R = 1$. If  $r(n,k)$ is the 
radius of the largest
ball inscribed in $\Delta_k^n = \Delta_k^n(1)$, then 
$r(n,\cdot) : (-\infty,
(\pi/2)^2) \to (0,1)$  is a strictly increasing continuous 
function for each $n
\ge 2$. Let $k(n)$ be determined by  $r(n,k(n)) = 1/2$.  
The optimality
question of $(*) $ is then resolved for $q = n+1$ 
according to the following. 

\thm {Theorem B}  Fix an integer $n \ge 2$ and $k < 
(\frac\pi{4})^2$.  For any
closed Riemannian $n$-manifold $M$ with $\Sec M \ge k$ and 
$\rad M \le 1$, 
$$
	\xt_{n+1} M < \xt_{n+1} D_k^n(1), 
$$
and this inequality is optimal if and only if $k \le k(n)$.

	There is an $\varepsilon(n) > 0$, so that if $k = k(n)$
and
$$
	\xt_{n+1} M \ge \xt_{n+1} D_{k(n)}^n(1) - \varepsilon(n),
$$
then $M$ is homeomorphic to $S^n$.  
\ethm

In contrast to the proof of Theorem A, this result is 
proved using  convergence
techniques.  It follows from these techniques via \cite 
{GPW} (cf.\ also a
recent result of Perelman announced in \cite {BGP}) that 
there are at most
finitely many topological types of manifolds $M$ as in 
Theorem B, for which
$\xt _{n+1} M \ge \text{xt}_{n+1} D_k^{n-1}(1) + 
\varepsilon$, for any fixed 
$\varepsilon > 0$.  

Once it has been observed that the double $D\Delta_k^n = 
\Delta_k^n
\coprod \Delta_k^n / (\partial \Delta_k^n \sim 
\partial\Delta_k^n)$ of the
simplex $\Delta_k^n$ has curvature $\curv D\Delta_k^n \ge 
k$ in distance
comparison sense and $\rad D\Delta_k^n = 1$ when $k \le 
k(n)$, the optimality
statement in Theorem B is fairly trivial.  The hard part 
is to prove the
recognition statement when $k = k(n)$ and the 
nonoptimality statement when $k
> k(n)$.  This, on the other  hand, follows from Theorem C 
below together with
results and tools developed in \cite {GP3, GPW}.

Following the terminology of \cite {BGP}, a FSCBB, or an 
{\it
Aleksandrov space\/} is a complete inner metric space with 
finite Hausdorff
dimension,  which is curved from below in (local) distance 
comparison sense.

\thm	{Theorem C}  Let $X$ be an $n$-dimensional 
Aleksandrov space with
$\curv X \ge k$, $k < (\frac\pi{4})^2, \rad X = 1$, and 
$\xt_{n+1} X =
\xt_{n+1} D_k^n(1)$.  Then
\roster
\item"{(i)}"  There is an isometric embedding of 
$\Delta_k^n$ in $X$ with
totally geodesic interior\/{\rm ;}
\item"{(ii)}"   If $X$ is a Poincar\'e duality space then 
$k \le k(n)$\/{\rm ;}
\item"{(iii)}"  If $X$ is a Poincar\'e duality
space and $k = k(n)$, then $X$ is isometric to 
$D\Delta_{k(n)}^n$.
\endroster
\ethm

The essential new technical tool used in the proof of 
Theorem C is the
following analogue of the rigidity version of Toponogov's 
distance comparison 
theorem for Aleksandrov spaces.

\thm{Theorem D}  Let $X$ be an Aleksandrov space with 
$\curv X \ge k$.  For any
pair $(p_0, c_0)$, where $c_o$ is a minimal geodesic in 
$X$ with end points
$p_1,p_2$ and $p_0 \notin c_0$, let $(\bar p_0, \bar c_0)$ 
be the
corresponding  pair in $S_k^2$, i.e.\  $\dist(p_i,p_j) = 
\dist(\bar p_i, \bar
p_j), 0 \le i < j \le 2$ {\rm (}if $k > 0$ assume all 
distances $< \pi/\sqrt
k)$.  Then for corresponding interior points $q \in c_0, 
\;  \bar q \in \bar
c_0$, we have 

{\rm (i)} $\dist(p_0,q) \ge d(\bar p_0, \bar q)$, and

{\rm (ii)} if equality holds, any minimal geodesic $c_q$ 
from $p_0$ to
$q$ spans together with $c_0$ a unique triangular surface 
isometric to the one
spanned by $(\bar p_0, \bar c_0)$ in $S_k^2$, and whose 
interior is totally
geodesic. 
\ethm

	Part (i) of this global distance comparison theorem was 
proved in \cite
{BGP}.  A rigidity statement for hinges as in Riemannian 
geometry (cf.\ \cite
{CE, G1}) follows by applying (ii) above twice.

	In the generality of Theorem C as stated above, we also 
need a new
metric characterization of $S_1^n$.  In order to describe 
this, let the
$q$-{\it packing  radius}, $\pack _qX$, of $X$ be the 
largest $r$ so that there
are $q$ disjoint open $r$-balls in $X$, i.e., 
$$
	2 \pack _qX = \max_{(x_1,\dotsc,x_q)} \; \min_{i<j} \; 
\text{dist}
(x_i,x_j),  \qquad  \; x_i,\dotsc,x_q \in X.
$$

\thm{Theorem E}  Let $X$ be an $n$-dimensional Aleksandrov 
space with $\curv 
(X) \ge 1$. Then
$$
	\pack _{n+2}(X) \le \pack_{n+2}(S_1^n)
$$
and equality holds if and only if $X$ is isometric to 
$S_1^n$.
\ethm

	This characterization together with Yamaguchi's fibration 
theorem 
\cite {Y}  also yields the following {\it pinching 
theorem\/} of independent
interest. 

\thm{Corollary F}  For each integer $n \ge 2$ there is an 
$\varepsilon =
\varepsilon(n) > 0$ such that any Riemannian $n$-manifold 
$M$ with $\Sec M \ge
1$ and $\pack_{n+2}(M)  \ge \pack_{n+2}(S_1^n) - 
\varepsilon$ is
diffeomorphic to $S^n$. \ethm

We conclude with another pinching---or rather {\it 
recognition
theorem\/}.  As in our previous results above, this is 
also based on a complete
metric classification of a certain class of Aleksandrov 
spaces.  In this case
the class consists of all $n$-dimensional Aleksandrov 
spaces $X$ with $\curv
(X) \ge 1$ and $\diam (X) = \xt_{n+1}(X) = \pi/2$.  Rather 
than giving the list
here, we point out that for each $n$ only two of them are 
Poinc{\'a}re duality
spaces, namely, $\Bbb RP_1^n = S_1^n /\Bbb Z_2$ and the  
double spherical
simplex $D\Delta_1^n$.  From this one derives

\thm {Theorem G}  For each $n \ge 2$ there is an 
$\varepsilon = \varepsilon(n)
> 0$ with the following property.  Any closed Riemannian 
$n$-manifold $M$  with
$\Sec M \ge 1$ and $\xt_{n+1}(M) \ge \pi/2 - \varepsilon$ 
is either
diffeomorphic to $\Bbb RP^n$ or  homeomorphic to $S^n$.  
Moreover, there are
metrics on $M \sim S^n, \Bbb RP^n$ with $\Sec M \ge 1$ , 
$\diam M \le \pi/2$, 
and  $\xt_{n+1}(M)$ arbitrarily close to $\pi/2$. \ethm

	This result can be viewed as a generalization of the main 
results in
\cite {GP1, OSY}.

	Details and further applications will be published in 
\cite {GM}.

	It is our pleasure to thank S. Ferry for helpful 
suggestions related to
the proof of Theorem C.

\enddocument